\documentclass[b5paper,oneside,10pt]{article}

\makeatletter

\newif\if@final
\newif\if@finalweb
\newif\if@proofs
\newif\if@proofsline
\newif\if@oldstyle
\newif\if@imfm
\@finalfalse
\@finalwebfalse
\@proofsfalse
\@imfmfalse

\if@proofsline\RequirePackage[displaymath, mathlines]{lineno}\fi
\RequirePackage{amsmath}
\RequirePackage{amsthm}
\RequirePackage{amsfonts}
\RequirePackage{geometry}
\if@final\RequirePackage{graphicx}\RequirePackage{lastpage}\RequirePackage{ccicons}\fi
\if@proofs\RequirePackage{eso-pic}\fi
\if@finalweb
\if@proofs
\RequirePackage[hyperfootnotes=false,colorlinks=true,pagebackref=true]{hyperref}
\else
\RequirePackage[hyperfootnotes=false,colorlinks=true]{hyperref}
\fi
\else
\RequirePackage[hyperfootnotes=false,colorlinks=true]{hyperref}
\fi

\numberwithin{equation}{section}

\RequirePackage{times}

\newcount\@authorcount
\newcount\@footcount
\if@oldstyle
\else
\fi
\if@final\else\fi
\def\@emails{}
\def\@authors{}
\def\p@pertitle{}
\def\blfootnote{\xdef\@thefnmark{}\@footnotetext} 
\long\def\symbolfootnote[#1]#2{\begingroup%
\def\thefootnote{\fnsymbol{footnote}}\footnote[#1]{#2}\endgroup} 
\def\@startpage{1}
\def\@endpage{x}
\def\@oddrunninghead{ }
\def\@evenrunninghead{ }
\def\@onlinedate{*}
\def\@receiveddate{*}
\def\@acceptdate{*}

\newcommand{\journalnumber}[1]{\def\@journalnumber{#1}}
\newcommand{\journalyear}[1]{\def\@journalyear{#1}}
\newcommand{\journalpages}[2]{%
\def\@startpage{#1}
\def\@endpage{#2}
}
\newcommand{\oddrunninghead}[1]{\def\@oddrunninghead{#1}}
\newcommand{\evenrunninghead}[1]{\def\@evenrunninghead{#1}}
\newcommand{\dates}[3]{%
\def\@receivedate{#1}%
\def\@acceptdate{#2}%
\def\@onlinedate{#3}}

\newenvironment{frontmatter}
{\thispagestyle{papertitle}%
\setcounter{page}{\@startpage}%
\vskip 10pt
\if@proofsline\global\linenumbers\fi%
}
{\vskip 20pt%
\blfootnote{\raggedright \ifnum\@authorcount=1\textit{E-mail addresses:}\else\textit{E-mail addresses:}\fi ~\@emails}%
\if@proofsline\global\linenumbers\fi%
}

\newcommand{\footadvance}[1]{%
\global\advance\@footcount by #1
}

\newcommand{\titledata}[2]{%
\global\def\p@pertitle{#1}%
\vskip 5mm
\begin{center}\LARGE\bfseries%
#1\if!#2!\else
\global\advance\@footcount by 1%
\symbolfootnote[\@footcount]{#2}\fi\end{center}%
\vskip 1mm
}

\newcommand{\authordata}[4]{
\advance\@authorcount by 1%
\xdef\@emails{\@emails\ifnum\@authorcount>1 , \fi {#3} (#1)}%
\xdef\@authors{\@authors\ifnum\@authorcount>1 , \fi {#1}}
\begin{center}%
\if!#4!
{\large #1}\\[1mm]
{\itshape #2}
\else
\global\advance\@footcount by 1%
{\large #1} \symbolfootnote[\@footcount]{#4}\\[1mm]
{\itshape #2}%
\fi
\end{center}}

\newcommand{\authordatanoemail}[3]{
\begin{center}%
\if!#3!
{\large #1}\\[1mm]
{\itshape #2}
\else
\global\advance\@footcount by 1%
{\large #1} \symbolfootnote[\@footcount]{#3}\\[1mm]
{\itshape #2}%
\fi
\end{center}}

\newcommand{\authordataff}[5]{
\advance\@authorcount by 1%
\xdef\@emails{\@emails\ifnum\@authorcount>1 , \fi #3 (#1)}%
\xdef\@authors{\@authors\ifnum\@authorcount>1 , \fi {#1}}
\begin{center}%
\global\advance\@footcount by 1%
{\large #1} \symbolfootnote[\@footcount]{#4}
\global\advance\@footcount by 1%
\symbolfootnote[\@footcount]{#5}
\\[1mm]
{\itshape #2}%
\end{center}}

\newcommand{\authordatawithfootnote}[4]{
\advance\@authorcount by 1%
\xdef\@emails{\@emails\ifnum\@authorcount>1 , \fi #3 (#1)}%
\xdef\@authors{\@authors\ifnum\@authorcount>1 , \fi {#1}}
\begin{center}%
{\large #1} $^{#4}$\\[1mm]
{\itshape #2}%
\end{center}}

\newcommand{\authordatanoaffil}[3]{
\advance\@authorcount by 1%
\xdef\@emails{\@emails\ifnum\@authorcount>1 , \fi #2 (#1)}%
\xdef\@authors{\@authors\ifnum\@authorcount>1 , \fi {#1}}
\if!#3!
{\large #1}
\else
\global\advance\@footcount by 1%
{\large #1} \symbolfootnote[\@footcount]{#3}%
\fi}

\newcommand{\authordatanoaffilnofootnote}[2]{
\advance\@authorcount by 1%
\xdef\@emails{\@emails\ifnum\@authorcount>1 , \fi #2 (#1)}%
\xdef\@authors{\@authors\ifnum\@authorcount>1 , \fi {#1}}
{\large #1\!}
}

\newcommand{\keywords}[1]{\def\@keywords{#1}}
\newcommand{\msc}[1]{\def\@msc{#1}}

\renewenvironment{abstract}
{
\hrule height 0.25pt
\vskip 5pt
\noindent \textbf{Abstract}
\vskip 5pt
}
{
\vskip 5pt
\noindent \textit{\small Keywords:~\@keywords}

\vskip 3pt
\noindent \textit{\small Math.\ Subj.\ Class.: \@msc}
\vskip 5pt
\hrule height 0.25pt
} 

\renewcommand{\section}{%
\@startsection{section}{1}{0pt}{-\baselineskip}{0.5\baselineskip}{\large\bfseries}%
}

\renewcommand{\subsection}{%
\@startsection{subsection}{2}{0pt}{-\baselineskip}{0.5\baselineskip}{\bfseries}%
}

\theoremstyle{plain}
\newtheorem{thm}{Theorem}[section]
\newtheorem{theorem}[thm]{Theorem}

\newtheorem{lemma}[thm]{Lemma}

\newtheorem{corollary}[thm]{Corollary}

\theoremstyle{definition}

\if@proofs
  \newcommand{\ps@pagestyle}{%
    \renewcommand\@oddhead{\hfil \small \textit{\@oddrunninghead} \hfil \textit{\thepage}}
    \renewcommand\@evenhead{\small \textit{\thepage} \hfil \textit{\@evenrunninghead} \hfil}
    \renewcommand\@evenfoot{}
    \renewcommand\@oddfoot{}
  }
  \newcommand{\ps@papertitle}{%
    \renewcommand\@oddhead{}
    \let\@evenhead\@oddhead
    \let\@oddfoot\@evenfoot
  }
\else
  \newcommand{\ps@pagestyle}{%
    \renewcommand\@oddhead{\hfil \small \textit{\@oddrunninghead} \hfil \textit{\thepage}}
    \renewcommand\@evenhead{\small \textit{\thepage} \hfil \textit{\@evenrunninghead} \hfil}
    \renewcommand\@evenfoot{}
    \renewcommand\@oddfoot{}
  }
  \newcommand{\ps@papertitle}{%
    \renewcommand\@oddhead{}
    \let\@evenhead\@oddhead
    \let\@oddfoot\@evenfoot
  }
\fi

\def\enumerate{%
  \ifnum \@enumdepth >\thr@@\@toodeep\else
    \advance\@enumdepth\@ne
    \edef\@enumctr{enum\romannumeral\the\@enumdepth}%
      \expandafter
      \list
        \csname label\@enumctr\endcsname
        {\usecounter\@enumctr\def\makelabel##1{\hss\llap{##1}}%
          \itemsep=0pt}%
  \fi}
\def\itemize{%
  \ifnum \@itemdepth >\thr@@\@toodeep\else
    \advance\@itemdepth\@ne
    \edef\@itemitem{labelitem\romannumeral\the\@itemdepth}%
    \expandafter
    \list
      \csname\@itemitem\endcsname
      {\def\makelabel##1{\hss\llap{##1}}%
        \itemsep=0pt}%
  \fi}

\renewenvironment{thebibliography}[1]
     {\section*{\refname}%
      \list{\@biblabel{\@arabic\c@enumiv}}%
           {\settowidth\labelwidth{\@biblabel{#1}}%
            \leftmargin\labelwidth
            \advance\leftmargin\labelsep
            \@openbib@code
            \usecounter{enumiv}%
            \let\p@enumiv\@empty
            \renewcommand\theenumiv{\@arabic\c@enumiv}
            \itemsep=0.0pt}%
      \sloppy
      \small
      \clubpenalty4000
      \@clubpenalty \clubpenalty
      \widowpenalty4000%
      \sfcode`\.\@m}
     {\def\@noitemerr
       {\@latex@warning{Empty `thebibliography' environment}}%
      \endlist}
\let\@openbib@code\@empty

\if@proofs
\newcommand{\BackgroundPic}{
\put(0,-50){
\parbox[b][\paperheight]{\paperwidth}{%
\vfill
\centering
\includegraphics[width=\paperwidth,height=\paperheight,keepaspectratio]{proofs.pdf}%
\vfill
}}}
\fi

\newcommand{\pagestyleend}
{
\if@final
\newwrite\contentsfile
\immediate\openout\contentsfile=contentsfile.txt
\makeatletter
\immediate\write\contentsfile{\@backslashchar contentsline\@charlb subsection\@charrb
\@charlb\@backslashchar textbf\@charlb\p@pertitle\@charrb\@backslashchar\@backslashchar}
\immediate\write\contentsfile{\@authors\@charrb\@charlb\@startpage\@charrb}
\makeatother
\closeout\contentsfile
}

\if@proofs\AddToShipoutPicture{\BackgroundPic}\fi

\newcommand{\fourIdx}[5]{%
\ensuremath{\mathchoice%
{{\setbox1=\hbox{$\scriptstyle{#1}$}%
\setbox2=\hbox{$\scriptstyle{#2}$}%
\setbox5=\hbox{$\displaystyle{#5}$}%
\hspace*{\ifnum\wd1>\wd2\wd1\else\wd2\fi}%
\copy5^{\scriptstyle\hspace{-\wd1}\hspace{-\wd5}#1\hspace{\wd5}#3}%
_{\scriptstyle\hspace{-\wd2}\hspace{-\wd5}#2\hspace{\wd5}#4}%
}}
{{\setbox1=\hbox{$\scriptstyle{#1}$}%
\setbox2=\hbox{$\scriptstyle{#2}$}%
\setbox5=\hbox{$\textstyle{#5}$}%
\hspace*{\ifnum\wd1>\wd2\wd1\else\wd2\fi}%
\copy5^{\scriptstyle\hspace{-\wd1}\hspace{-\wd5}#1\hspace{\wd5}#3}%
_{\scriptstyle\hspace{-\wd2}\hspace{-\wd5}#2\hspace{\wd5}#4}%
}}
{{\setbox1=\hbox{$\scriptscriptstyle{#1}$}%
\setbox2=\hbox{$\scriptscriptstyle{#2}$}%
\setbox5=\hbox{$\scriptstyle{#5}$}%
\hspace*{\ifnum\wd1>\wd2\wd1\else\wd2\fi}%
\copy5^{\scriptscriptstyle\hspace{-\wd1}\hspace{-\wd5}#1\hspace{\wd5}#3}%
_{\scriptscriptstyle\hspace{-\wd2}\hspace{-\wd5}#2\hspace{\wd5}#4}%
}}
{{\setbox1=\hbox{$\scriptscriptstyle{#1}$}%
\setbox2=\hbox{$\scriptscriptstyle{#2}$}%
\setbox5=\hbox{$\scriptscriptstyle{#5}$}%
\hspace*{\ifnum\wd1>\wd2\wd1\else\wd2\fi}%
\copy5^{\scriptscriptstyle\hspace{-\wd1}\hspace{-\wd5}#1\hspace{\wd5}#3}%
_{\scriptscriptstyle\hspace{-\wd2}\hspace{-\wd5}#2\hspace{\wd5}#4}%
}}
}}

\usepackage{amssymb}
\usepackage[displaymath, mathlines]{lineno}
\usepackage[dvipsnames]{xcolor}
\usepackage{graphicx, xfrac, lineno, float, subcaption, tasks, comment, booktabs, multirow}
\usepackage[normalem]{ulem}
\usepackage{datetime}
\usepackage[super]{nth}
\usepackage{colortbl}
\usepackage{enumerate}

\makeatother

\usepackage{mathtools}
\usepackage{tikz}
\usepackage{enumitem}

\definecolor{orchid}{rgb}{0.85, 0.44, 0.84}

\newcommand{\Hclass}{\mathcal H}
\newcommand{\TT}{\mathcal T}
\newcommand{\RR}{\mathcal R}
\newtheorem{question}[thm]{Question}
\DeclareMathOperator{\perm}{per}
\newcommand{\hnum}{\chi_h^{(2)}}
\begin{document}

\begin{frontmatter}

\titledata{On 2-Homogeneous Colourings of Cubic Bipartite Graphs}{}

\authordata{Rong-Xia Hao}{School of Mathematics and Statistics, Beijing Jiaotong University, Beijing, 100044, PR China}{rxhao@bjtu.edu.cn}{}
\authordata{Jianbing Liu}{School of Mathematics and Statistics, Beijing Jiaotong University, Beijing, 100044, PR China}{jbliu1@bjtu.edu.cn}{}
\authordata{Federico Romaniello}
{Dipartimento per l'Innovazione Umanistica, Scientifica e Sociale\\ Universit\`{a} degli Studi della Basilicata, Potenza, Italy\\ }{federico.romaniello@unibas.it}{}

\msc{05C15, 05C70}

\begin{abstract}
A proper vertex colouring is 2-homogeneous if every open neighbourhood contains exactly two colours. We show that every cubic bipartite graph admits a 2-homogeneous colouring with at most four colours, proving a conjecture of Janicová, Madaras, Soták and Lu\v{z}ar. For 3-connected det-extremal cubic bipartite graphs, we establish the stronger bound $\chi_h^{(2)}(G) = 3$, showing that every such graph admits a fall 3-colouring.
\end{abstract}
\end{frontmatter}

\section{Introduction}\label{sec:intro}

All graphs considered are finite, simple and connected unless stated
otherwise. Connectivity means vertex connectivity. A proper vertex
colouring $c$ of a graph $G$ is \emph{2-homogeneous} if
\[
 |c(N_G(v))|=2\qquad\text{for every }v\in V(G),
\]
where $N_G(v)$ denotes the open neighbourhood of $v$. When such a
colouring exists, let $\hnum(G)$ be the minimum number of colours it
uses. Every such colouring uses at least three colours. For a cubic graph, the colours on the
three neighbours of each vertex have multiplicities two and one.

Janicov\'a, Madaras, Sot\'ak and Lu\v{z}ar proved that every cubic
bipartite graph admits a 2-homogeneous colouring with at most six
colours and conjectured that four colours suffice
\cite[Theorem~1.2 and Conjecture~4.1]{JMSL}. They also proved the
four-colour bound for regular bipartite graphs of degree at least four.
Our main result settles their conjecture and answers Problem~5.1 of
Lu\v{z}ar and Sot\'ak~\cite{CC} for cubic bipartite graphs.

\begin{theorem}\label{main}
Every cubic bipartite graph $G$ satisfies $3\le\hnum(G)\le4$.
The upper bound is attained by $K_{3,3}$.
\end{theorem}

The sharpness example is already observed in~\cite[p.~358]{JMSL}.
Each part of $K_{3,3}$ is the neighbourhood of every vertex in the
other part and must therefore use exactly two colours. Properness
forces these two pairs of colours to be disjoint. Conversely,
assigning two colours to each part gives a 2-homogeneous four-colouring.

Let $B$ be the biadjacency matrix of a cubic bipartite graph. Its two
parts have the same size. The graph is \emph{det-extremal} if
$|\det B|=\perm B$, where $\perm B$ denotes the permanent of $B$.
Equivalently, for fixed orders of the two parts, all perfect matchings
contribute with the same sign to $\det B$.

A proper 2-homogeneous colouring with three colours is precisely a
\emph{fall 3-colouring}: every closed neighbourhood contains all three
colours~\cite{Dunbar}. Equivalently, its three colour classes are
independent dominating sets. Our proof gives the following
three-colour result.

\begin{corollary}\label{cor:three}
Let $G$ be a 3-connected det-extremal cubic bipartite graph.
Then $G$ admits a fall 3-colouring. In particular, $\hnum(G)=3$.
\end{corollary}

The non-det-extremal case of Theorem~\ref{main} is a consequence of
Choi, Kwon and Park's characterisation in terms of the total domatic
number~\cite[Theorem~1.4]{CKP}. Their theorem provides a partition into
two total dominating sets. Intersecting these sets with the bipartition
gives the four colour classes. For det-extremal graphs of connectivity
three, we use McCuaig's construction from the Heawood graph and prove
an extension lemma with prescribed boundary colours. The case of
connectivity two follows from a cut lemma of Labbate and
Romaniello~\cite{LR}, based on the structure theorem of Funk et
al.~\cite{FJLS}.

Section~\ref{sec:nonext} gives the non-det-extremal construction.
Section~\ref{sec:proof} treats the det-extremal case and completes the
proofs of Theorem~\ref{main} and Corollary~\ref{cor:three}; the boundary
extension lemma is proved in Section~\ref{sec:boundaryproof}.
Section~\ref{sec:conclusions} gives a sharp example of connectivity two
and relates the result to the question for general cubic graphs.

\section{Non-det-extremal graphs}\label{sec:nonext}

A set $D\subseteq V(G)$ is \emph{total dominating} if every vertex
of $G$ has a neighbour in $D$, including the vertices of $D$ itself.
The \emph{total domatic number} $d_t(G)$ is the maximum number of
pairwise disjoint total dominating sets. Choi, Kwon and
Park~\cite[Theorem~1.4]{CKP} prove that a connected cubic graph is
det-extremal if and only if its total domatic number is one.
For bipartite graphs, their definition using the adjacency matrix
is equivalent to the definition using the biadjacency matrix above.
We use the following consequence.

\begin{lemma}\label{nonext}
Every non-det-extremal cubic bipartite graph admits a
proper 2-homogeneous colouring with at most four colours.
\end{lemma}

\begin{proof}
Let $V(G)=X\sqcup Y$ be the bipartition. Since $G$ is not
det-extremal, $d_t(G)\ge2$. Choose disjoint total dominating sets
$D_0,D_1$ and add any remaining vertices to $D_1$. Both sets remain
total dominating and now partition $V(G)$. Define
\[
c(v)=
\begin{cases}
1,&v\in X\cap D_0,\\
2,&v\in X\cap D_1,\\
3,&v\in Y\cap D_0,\\
4,&v\in Y\cap D_1.
\end{cases}
\]
Every edge joins $X$ to $Y$, whose palettes are disjoint, so the
colouring is proper. Every vertex has a neighbour in each of
$D_0,D_1$. Hence each vertex of $X$ sees exactly the colours $3,4$,
and each vertex of $Y$ sees exactly the colours $1,2$.
\end{proof}

Equivalently, one may apply the incidence criterion
of~\cite[Theorem~3.5]{CKP} to the neighbourhood hypergraph on each
part. A direct construction using perfect matchings and an even
directed cycle is given in~\cite[Lemmas~3.3--3.4]{LR}.

\section{Det-extremal graphs}\label{sec:proof}

To prove Theorem~\ref{main} and Corollary~\ref{cor:three}, we first
construct three-colourings for the 3-connected det-extremal graphs.
We then use a 2-edge cut to obtain the four-colour bound in
connectivity two.

\subsection{Vertex-sums and boundary states}

Let $H$ be the Heawood graph, the incidence graph of the Fano plane.
Each of its seven lines contains three points, any two points determine
a unique line, and any two lines meet in a unique point. A
\emph{vertex-sum} deletes one vertex from each of two disjoint cubic
graphs and joins their two sets of former neighbours by a matching.
Let $\Hclass$ denote the class generated from $H$ by vertex-sums.
By McCuaig's theorem~\cite{McC}, the graphs in $\Hclass$ are precisely
the 3-connected det-extremal cubic bipartite graphs
\cite[Theorem~2.2]{FJLS}.

It therefore suffices to show that every graph in $\Hclass$ admits
a fall 3-colouring. We shall prove a boundary extension lemma that
allows these colourings to be combined across vertex-sums.

For the Heawood construction, we use the colour set $\{0,1,2\}$. Deleting a vertex
from a cubic graph and retaining its three incident half-edges gives a
\emph{3-pole}. Its three boundary vertices are called \emph{ports}. Fix an
ordering of the half-edges. A \emph{boundary state} is a pair $(E,P)$ of
triples, where $p_i$ is the colour at the $i$th port and $e_i$ is the prescribed
colour at the other end of its half-edge. We require $e_i\ne p_i$ for every
$i$. The 3-pole \emph{realises} $(E,P)$ if it has a proper colouring with
these port colours in which every vertex sees both other colours, counting
$e_i$ at the $i$th port.

Let $q(W)$ denote the number of distinct entries of a triple $W$. Among the
proper boundary states, define
\begin{align*}
 (E,P)\in\TT &\quad\Longleftrightarrow\quad q(P)=1,\quad q(E)=2,\\
 (E,P)\in\RR &\quad\Longleftrightarrow\quad q(P)\ge2,\quad
                         (q(P),q(E))\ne(3,3).
\end{align*}
A state in $\TT$ describes a single coloured vertex: its own colour is
repeated at its three incident ends, and the opposite ends use the other two
colours. The states in $\RR$ will be realised by the 3-poles used in the
construction. In particular, if $q(P)=2$, every proper choice of $E$ is
allowed. The state $(000,112)$ permits the deleted vertex to be restored
with colour $0$.

We shall also use colours on individual edge ends. Let $F=H-z$ and
$S\subseteq V(F)$. Assign a colour to each end of each internal edge and
prescribe the two colours at every boundary half-edge. Opposite ends must
have distinct colours. At $u\in V(F)$, let $P_u$ be the triple of colours
at its own ends and let $E_u$ be the corresponding triple at the opposite
ends. An \emph{assignment for $S$} satisfies
\[
 (E_u,P_u)\in
 \begin{cases}
  \RR,&u\in S,\\
  \TT,&u\notin S.
 \end{cases}
\]
Thus every vertex outside $S$ has a single colour, whereas a vertex of $S$
may be replaced by a 3-pole realising its assigned state. At this stage,
the three ends at a vertex of $S$ need not have the same colour.

\subsection{The boundary extension lemma}\label{sec:boundaryproof}

We first show how an assignment for an independent set extends to
an arbitrary selection of vertices. We then study recolourings of
individual edge ends in a fall 3-colouring of $F$. An auxiliary
multigraph determines when all selected vertices can be recoloured.
Finally, explicit initial colourings cover the seven types of
boundary state.

\subsubsection{Reduction to an independent set}
\begin{lemma}\label{lem:activation}
Let $S\subseteq V(F)$, and let $I$ be the set of isolated vertices
of $F[S]$. Every assignment for $I$ can be changed into an assignment
for $S$ with the same boundary, by recolouring only ends of edges
of $F[S\setminus I]$.
\end{lemma}
\begin{proof}
Put $U=S\setminus I$ and take an assignment for $I$ with the
prescribed boundary. For each $u\in U$, let $c_u$ be the
common colour at its three incident ends.

Choose a maximal matching $M$ of $F[U]$. On each edge $uv\in M$,
interchange the colours at its ends. The edge remains proper, and
each vertex covered by $M$ now has two colours at its own ends.
Since $F[U]$ has no isolated vertices and $M$ is maximal, every
uncovered vertex $u$ has a covered neighbour $w$. On $uw$, recolour
only the end at $u$ with the third colour distinct from $c_u,c_w$.
The end at $w$ still has colour $c_w$, since $uw\notin M$.

Now $q(P_u)=2$ for every $u\in U$. All the new states at these vertices
belong to $\RR$. Only ends of edges inside $F[U]$ were changed, so the
states at vertices outside $U$ and the prescribed boundary are preserved.
\end{proof}

\subsubsection{Activations}

Let $c$ be a fall 3-colouring of $F$ with prescribed boundary, again
counting the exterior colours at its ports. Let $S$ be independent.
On an internal edge $uv$ with $u\in S$ we may change the end at $u$
to the third colour distinct from $c(u)$ and $c(v)$. We call this
change an \emph{activation} at $u$. The end at $v$ retains colour
$c(v)$, and boundary half-edges are kept fixed.

For $u\in S$, the triple $E_u$ still uses two colours since the
neighbours of $u$ lie outside $S$. One or two activations at $u$ leave
both the original colour and a new colour in $P_u$. If all three ends
are changed, the two colours in $E_u$ give two different new colours
in $P_u$. Thus any nonempty set of activations at $u$ gives a state
in $\RR$. It remains to determine which changes preserve $\TT$ at the
vertices outside $S$.

For $v\notin S$, write its neighbour colours as $a,a,b$. Call the
neighbour of colour $b$ the \emph{minority neighbour} of $v$ and the
other two its \emph{majority neighbours}. An exterior neighbour can
have either role and is always outside $S$. Activating an incidence with $v$ changes its neighbour colour
from $a$ to $b$ or from $b$ to $a$. The
triple becomes monochromatic in just two cases, when only the minority
position is changed or when both majority positions and no other
position are changed.

Define a multigraph $D$ on $S$. For each $v\notin S$ with
$N(v)\cap S\ne\varnothing$, consider the following possibilities.
\begin{enumerate}[label=(\alph*)]
\item If $N(v)\cap S$ consists of the minority neighbour alone, its
incidence with $v$ cannot be activated.
\item If $N(v)\cap S$ consists of the two majority neighbours, join
them by an edge of $D$. At most one of the two incidences with $v$ may
be activated, and either choice is possible.
\item In every other case all incidences from $S$ to $v$ may be
activated together. Call each vertex of $N(v)\cap S$ \emph{free}.
\end{enumerate}

\begin{lemma}\label{lem:activation-criterion}
Every vertex of $S$ can be activated while preserving the boundary
state and the neighbourhood condition outside $S$ if and only if
every component of $D$ contains a free vertex or a cycle.
Two parallel edges count as a cycle.
\end{lemma}
\begin{proof}
Perform all the activations in (c). In a component containing a free
vertex, root a spanning tree at that vertex and direct its edges away
from the root. For each directed edge use (b) to activate its head.
In a component containing a cycle, extend the cycle to a connected
spanning subgraph with just that cycle. Direct the cycle cyclically
and all other edges away from it. Again every vertex requiring an
activation receives one.

Conversely, let $C$ be a tree component with no free vertex. Each
activation available in $C$ uses an edge of $C$ as in (b), and each
such edge can serve only one end. Since $|E(C)|=|V(C)|-1$, at least
one vertex cannot be activated.
\end{proof}

\subsubsection{Point and line obstructions}

By vertex transitivity of $H$ we may take $z$ to be the line $123$.
Label the six retained lines by
\[
a=145,\quad b=167,\quad c=246,\quad d=257,\quad e=347,\quad f=356.
\]
The points $1,2,3$ are the ports. The points $4,5,6,7$, which lie
outside the deleted line, are called the \emph{outside points}. For an independent set $S$, write $X$ for its set of points
and $Y$ for its set of lines. No point of $X$ lies on a line of $Y$. The auxiliary multigraph
$D$ is the disjoint union of $D[X]$ and $D[Y]$, so point and line
activations can be chosen independently.

Fix a fall 3-colouring with prescribed boundary. The \emph{point
majority graph} $M_P$ has the seven points as vertices and joins
the two equally coloured points on each retained line. For the
\emph{line majority graph} $M_L$, take the six retained lines and
three fixed vertices representing the exterior neighbours of the ports.
At each point, join its two equally coloured neighbours. A retained
line is \emph{marked} if it is joined in $M_L$ to a fixed vertex.
Whenever a marked line is selected, it is free by (c).

\begin{lemma}\label{lem:majority-obstructions}
For the fixed colouring, the vertices of $X$ cannot all be activated
precisely when one of the following holds.
\begin{enumerate}[label=(P\arabic*)]
\item $X$ is the vertex set of a tree component of $M_P$ containing
an outside point.
\item $X\subseteq\{1,2,3\}$ and $X$ contains an isolated vertex
of $M_P$.
\end{enumerate}
The vertices of $Y$ cannot all be activated precisely when $Y$ is an
entire colour class of the retained lines of size one or two and none
of its lines is marked.
\end{lemma}
\begin{proof}
First consider the point vertices. Let $C$ be a component of $D$ on
points containing no free vertex. If $xy$ is an edge of $M_P$ and
$x\in V(C)$, then $y$ must be selected and the minority point on
the line $xy$ must be unselected. Otherwise $x$ is free by (c).
Hence $xy\in E(D)$ by (b), and $y\in V(C)$. It follows that $C$
is a whole component of $M_P$ with the same edges.

Suppose $x\in V(C)$ is an outside point. Every other point lies on
a retained line with $x$. A point of the same colour is joined to $x$
in $M_P$. A selected point of a different colour would make $x$ free
by (c) on their common line. Hence $X=V(C)$ is the colour class of
$x$. Conversely, if $X$ is the vertex set of such a component of
$M_P$, its selected points on any line form the majority pair or
the minority singleton, whenever the line meets $X$. Thus $D[X]=C$
and has no free vertex. Lemma~\ref{lem:activation-criterion} gives (P1).

If $C$ contains only ports, it consists of a single isolated port
because two ports lie together only on the deleted line. Selecting
an outside point would make the port free on their common line.
Conversely, if only ports are selected, an isolated port has no
available activation. It is the minority point on each of its two
retained lines, and no other port lies on either. This proves (P2).

For the line vertices, any two retained lines meet in a point. Lines
of the same colour form the majority pair at that point, so every
line colour class induces a clique in $M_L$. As above, a component
of $D$ on lines with no free vertex contains the whole clique of its
colour. No line in it can be marked. A selected line of another colour
would give a free activation where it meets a line of the clique.
Thus $Y$ must be exactly one unmarked colour class. Conversely, for
such a class every majority pair in the clique gives an edge of $D$
and there are no free vertices. The clique is a tree exactly when it
has one or two vertices. The conclusion follows from
Lemma~\ref{lem:activation-criterion}.
\end{proof}

\subsubsection{Initial colourings of the Heawood 3-pole}

We may permute the three colours and the three ports. Indeed, the points
of the deleted line are the three nonzero vectors of a two-dimensional
subspace of $\mathbb F_2^3$. Every permutation of these points is linear
on the subspace and extends to $\mathbb F_2^3$ by fixing a vector
outside it. The resulting automorphism preserves incidence and hence
also preserves independence of the selected vertices.

There are seven types of boundary state up to these permutations.
If $P$ is bichromatic, we may take $P=112$. Up to interchanging the
first two ports, the first two entries of $E$ are $00$, $02$ or $22$,
and the last entry is $0$ or $1$. The rows $A_0,B,C,D,E_0,F_0$ of
Table~\ref{tab:seeds} give these six possibilities. With $P=112$, their
exterior triples are $000,221,001,020,220,021$, respectively. If $P$
has three colours, then $E$ has exactly two by properness and the
definition of $\RR$. Its repeated colour must occur at the two ports
whose own colours differ from it. Renaming the colours and ordering
the ports gives $(E,P)=(001,120)$.

Table~\ref{tab:seeds} gives the initial colourings for these cases.
The string $Q$ records the colours of the points $4,5,6,7$, and $L$
records the colours of the lines $a,b,c,d,e,f$, in this order.
In each row the points of a retained line have two colours and the
line has the third colour. Every point sees both colours different
from its own when the exterior neighbours of the ports are included.

For each colouring, the last two columns list the selections of points
and lines for which activation fails. We use
\[
 \mathcal P_* = \{X:\varnothing\ne X\subseteq\{1,2,3\}\},
 \qquad
 \mathcal P_2 = \{X:2\in X\subseteq\{1,2,3\}\}.
\]
The entries follow from Lemma~\ref{lem:majority-obstructions}. For
instance, in row $C$ the components of the point majority graph on
$\{4\}$ and $\{3,7\}$ are trees, whereas the component on
$\{1,2,5,6\}$ contains a cycle. The only unmarked line colour class
of size at most two is $\{a,c\}$. Thus row $C$ fails on the point
side exactly for $X=\{4\}$ or $X=\{3,7\}$, and on the line side
exactly for $Y=\{a,c\}$.

\begin{table}[htb]
\centering
\begin{tabular}{lcccc@{\qquad}ll}
\toprule
 & $E$ & $P$ & $Q$ & $L$ & Point obstructions & Line obstructions\\
\midrule
$A_0$ & $000$ & $112$ & $0000$ & $222211$ & $\mathcal P_*$ & $\{e,f\}$\\
$A_1$ & $000$ & $112$ & $0112$ & $202010$ & $\{4\},\ \{3,7\}$ & $\{a,c\},\ \{e\}$\\
$B$   & $001$ & $110$ & $0112$ & $202012$ & $\{7\},\ \{3,4\}$ & $\varnothing$\\
$C$   & $001$ & $112$ & $0112$ & $202010$ & $\{4\},\ \{3,7\}$ & $\{a,c\}$\\
$D$   & $001$ & $122$ & $0122$ & $201010$ & $\{4\},\ \{1,5\}$ & $\{a\}$\\
$E_0$ & $001$ & $220$ & $0221$ & $101021$ & $\{7\},\ \{3,4\}$ & $\{e\}$\\
$F_0$ & $012$ & $100$ & $0012$ & $202112$ & $\{7\},\ \{1,6\}$ & $\varnothing$\\
$J$   & $001$ & $120$ & $0101$ & $221022$ & $\mathcal P_2$ & $\{c\}$\\
\bottomrule
\end{tabular}
\caption{Initial colourings of the Heawood 3-pole. The last two
columns list the exact selected sets for which point or line
activation fails. An empty family means that every selection on
that side can be activated.}
\label{tab:seeds}
\end{table}

We shall also use the following automorphisms, each of which fixes the
three ports individually.
\[
\begin{array}{c|c|c}
 &\text{outside points}&\text{retained lines}\\\hline
\alpha&(4\ 5)(6\ 7)&(c\ d)(e\ f)\\
\beta &(4\ 6)(5\ 7)&(a\ b)(e\ f)\\
\gamma=\alpha\beta&(4\ 7)(5\ 6)&(a\ b)(c\ d)
\end{array}
\]
If $c$ is an initial colouring, then $\alpha c$ has the same boundary
state, and its obstructions are the images under $\alpha$ of those of
$c$. The same applies to $\beta$ and $\gamma$.

We are now ready to prove the boundary extension lemma.

\begin{lemma}\label{lem:finite}
For every $S\subseteq V(F)$ and every boundary state $(E,P)\in\RR$
of $F$, there is an assignment for $S$ with boundary $(E,P)$.
\end{lemma}

\begin{proof}
By Lemma~\ref{lem:activation}, it suffices to consider an independent
set $S$. Write $S=X\cup Y$, where $X$ consists of points and $Y$ of
lines. By Lemmas~\ref{lem:activation-criterion}
and~\ref{lem:majority-obstructions}, an initial colouring suffices whenever neither $X$ nor $Y$ is
one of its listed obstructions.

\emph{$(E,P)=(000,112)$.} Consider $A_0,A_1$ and $\alpha A_1$.
If $X\in\mathcal P_*$, then $X$ contains only ports, so it is a point
obstruction for neither $A_1$ nor $\alpha A_1$. Their line obstruction
families are $\{\{a,c\},\{e\}\}$ and $\{\{a,d\},\{f\}\}$, which
are disjoint. The other possible obstruction for $A_0$ is $Y=\{e,f\}$.
This set occurs in neither line column for $A_1$ and $\alpha A_1$,
and their point obstruction families $\{\{4\},\{3,7\}\}$ and
$\{\{5\},\{3,6\}\}$ are disjoint.

\emph{$(E,P)=(001,110)$ or $(012,100)$.} For the first state use
$B,\alpha B$, and for the second use $F_0,\alpha F_0$. In each pair
there is no line obstruction, and the two point obstruction families
are disjoint.

\emph{$(E,P)=(001,112)$.} The colourings $C$ and $\alpha C$ have
disjoint point obstruction families and disjoint line obstruction
families. For both colourings to fail, one of the following must hold.
\[
\begin{array}{ll}
 X\in\{\{4\},\{3,7\}\},&Y=\{a,d\},\\
 X\in\{\{5\},\{3,6\}\},&Y=\{a,c\}.
\end{array}
\]
Both alternatives contradict the independence of $S$. Indeed,
$4,5$ lie on $a$, $7$ lies on $d$, and $6$ lies on $c$.

\emph{$(E,P)=(001,122)$ or $(001,220)$.} For $(001,122)$, the
colourings $D,\alpha D$ have disjoint point obstruction families.
Their only common line obstruction is $Y=\{a\}$, so both can fail
only for this choice of $Y$. Similarly, $\beta D,\gamma D$ can both
fail only for $Y=\{b\}$. At least one of the four colourings therefore
suffices. For $(001,220)$ the same argument applies to the pairs
$E_0,\gamma E_0$ and $\alpha E_0,\beta E_0$. Their respective common
line obstructions are $\{e\}$ and $\{f\}$.

\emph{$(E,P)=(001,120)$.} If $2\notin S$, then $X$ is a point
obstruction for neither $J$ nor $\alpha J$. Their only line
obstructions are $\{c\}$ and $\{d\}$, respectively.
Suppose now that $2\in S$. Apply the case $(001,110)$ to
$S\setminus\{2\}$. The three ends at point $2$ then have colour $1$.
Recolour only the end at point $2$ of its boundary half-edge with
colour $2$. Its two internal ends retain colour $1$, and its exterior
neighbour retains colour $0$. The boundary half-edge remains proper,
while $q(P_2)=2$ gives a state in $\RR$. No internal end is
changed, and the new outer boundary state is $(001,120)$.

These cases cover every state in $\RR$ and complete the proof.
\end{proof}

\subsection{Extension through vertex-sums}

We now apply Lemma~\ref{lem:finite} to the vertex-sum construction
of $\Hclass$. Its prescribed states allow simultaneous substitutions,
including substitutions at adjacent vertices.

\begin{lemma}\label{lem:rooted}
Let $G\in\Hclass$ and let $r\in V(G)$. The 3-pole $G-r$ can be obtained
from vertex-deleted Heawood graphs by successive substitutions at
vertices. It realises every boundary state in $\RR$.
\end{lemma}
\begin{proof}
We first establish the recursive representation by induction on
the number of vertex-sums. The assertion is immediate for $G=H$.
Suppose that the final sum deletes $a\in V(A)$ and $b\in V(B)$,
and that $r$ survives on the $A$-side. Then $G-r$ is obtained from
$A-r$ by replacing $a$ with $B-b$. The induction hypothesis gives
both pieces with their attachments in the prescribed order. If $a$
is a port, its boundary half-edge is inherited by the corresponding
port of $B-b$.

We prove the assertion about boundary states by induction on this
representation. The base case $H-z$ is Lemma~\ref{lem:finite}
with $S=\varnothing$. In the induction step, consider the outer copy
of $H-z$ and let $S$ be the set of its vertices occupied by smaller
3-poles. Given a boundary state in $\RR$, Lemma~\ref{lem:finite}
assigns compatible endpoint colours with relation $\RR$ at the
vertices of $S$ and relation $\TT$ at the other vertices. By
induction, each smaller 3-pole admits the colouring specified at
its three ports. On every edge joining two pieces, the colour at
either end is the prescribed exterior colour for the other end.
The colourings therefore combine to give the required state on $G-r$.
\end{proof}

Restoring the deleted vertex gives the required three-colouring.

\begin{lemma}\label{lem:threeH}
Every graph in $\Hclass$ admits a fall 3-colouring.
\end{lemma}
\begin{proof}
Let $G\in\Hclass$ and choose $r\in V(G)$.
Lemma~\ref{lem:rooted} gives a colouring of $G-r$ with boundary
state $(E,P)=(000,112)$. Restore $r$ with colour $0$. Each port
then sees its prescribed exterior colour, and $r$ sees both $1$
and $2$. The resulting colouring is a fall 3-colouring of $G$.
\end{proof}

We can now prove Corollary~\ref{cor:three}.

\begin{proof}[Proof of Corollary~\ref{cor:three}]
By McCuaig's characterisation, $G\in\Hclass$, so
Lemma~\ref{lem:threeH} gives a fall 3-colouring. The lower bound
$\hnum(G)\ge3$ yields the equality.
\end{proof}

\subsection{Two-edge cuts and the four-colour bound}

For a 2-edge cut, the \emph{closure} of a shore is obtained by
joining its two boundary vertices with a new edge, called its
\emph{marker}. We use the following lemma of Labbate and
Romaniello~\cite[Lemma~3.18]{LR}. We include its proof, taking the
decomposition cactus of Funk et al.~\cite{FJLS} as a structural input.

\begin{lemma}\label{lem:terminal}
Let $G$ be a det-extremal cubic bipartite graph of connectivity two. It has a 2-edge cut with simple cubic bipartite closures $P,R$ such
that $P\in\Hclass$ and $R$ is not det-extremal.
\end{lemma}
\begin{proof}
For either shore, let $X',Y'$ be its bipartition classes and let
$d_X,d_Y$ count its outgoing edges incident with these classes.
Degree counting gives
$3(|X'|-|Y'|)=d_X-d_Y$, where $d_X+d_Y=2$.
Consequently $d_X=d_Y=1$: the boundary vertices lie in opposite
parts and the closure is cubic and bipartite. At this stage a
marker may be parallel to an original edge.

In the decomposition cactus, vertices represent classes joined
by three edge-disjoint paths. Closing a class with its markers
gives either a simple 3-connected cubic bipartite graph or a 3-bond, that
is, two vertices joined by three parallel edges. Choose a
non-cut vertex on an end cycle of the cactus. It has degree two,
and its closure $P$ has one marker. By
\cite[Theorem~4.10]{FJLS}, this piece is det-extremal.
It cannot be a 3-bond, since its two original edges would already
be parallel in $G$. Thus $P$ is simple and 3-connected, and
McCuaig's characterisation gives $P\in\Hclass$.

Let $R$ be the other closure, with markers $e=uv$ in $P$ and
$f=xy$ in $R$, so that the cut edges are $uy,vx$.
Choose perfect matchings $p_1,p_0$ of $P$ containing and avoiding
$e$, respectively, and similarly $r_1,r_0$ for $f$ in $R$.
Such choices exist by three-edge-colourability of cubic
bipartite multigraphs. Let $C_P,C_R$ be the cycles containing
the markers in $p_1\triangle p_0$ and $r_1\triangle r_0$.
Switching an alternating cycle of length $2\ell$ permutes the
matched images cyclically and changes the matching sign by
$(-1)^{\ell-1}$. Since $P$ is det-extremal,
$|C_P|\equiv2\pmod4$.

Replace the markers in $C_P\cup C_R$ by $uy,vx$. The resulting
cycle has length $|C_P|+|C_R|$ and alternates with respect to
\[
 M=((p_1\cup r_1)\setminus\{e,f\})\cup\{uy,vx\}.
\]
Det-extremality of $G$ therefore gives
$|C_P|+|C_R|\equiv2\pmod4$, whence $|C_R|\equiv0\pmod4$.

If $f$ were parallel to an original edge $g$, we could take
$r_0=(r_1\setminus\{f\})\cup\{g\}$. Then $C_R$ would have
length two, contradicting the congruence just obtained.
These are the only possible parallel edges in $R$, so $R$ is
simple. Finally, switching $C_R$ reverses a perfect-matching
sign, proving that $R$ is not det-extremal.
\end{proof}

We shall combine the colourings of the two closures using the
following elementary lemma.

\begin{lemma}\label{lem:glue}
Let $G_1,G_2$ be disjoint graphs with proper 2-homogeneous
colourings using at most $k$ colours. Choose edges $uv\in E(G_1)$
and $xy\in E(G_2)$. The graph obtained by deleting $uv,xy$ and
adding $uy,vx$ also has a proper 2-homogeneous colouring using at
most $k$ colours.
\end{lemma}
\begin{proof}
Use the same $k$-colour palette on both graphs. Since the ends of
each chosen edge have different colours, a permutation of the
palette on $G_1$ makes $c(u)=c(x)$ and $c(v)=c(y)$.
The two new edges are proper. At each of their four ends, the
colour of the deleted neighbour equals that of the new neighbour.
Every multiset of neighbour colours is therefore unchanged.
\end{proof}

We are now ready to prove Theorem~\ref{main}.

\begin{proof}[Proof of Theorem~\ref{main}]
The lower bound and equality for $K_{3,3}$ were established in
Section~\ref{sec:intro}. If $G$ is not det-extremal, apply
Lemma~\ref{nonext}.
A cubic bipartite graph has no bridge: counting degrees on one
shore of a bridge would give $3(|X'|-|Y'|)=\pm1$.
A cut vertex in a cubic graph would give a bridge, so $G$ has
connectivity two or three.
If $G$ is det-extremal and 3-connected, apply Lemma~\ref{lem:threeH}.
Otherwise choose the closures $P,R$ of Lemma~\ref{lem:terminal}.
Lemma~\ref{lem:threeH} colours $P$ with three colours, and
Lemma~\ref{nonext} colours $R$ with at most four colours.
Regard both as colourings from the same four-colour palette. Let $uv,xy$ be their markers,
labelled so that the original cut edges are $uy,vx$.
Lemma~\ref{lem:glue} combines the two colourings.
\end{proof}

\section{Concluding remarks}\label{sec:conclusions}

The four-colour bound remains sharp at connectivity exactly two. Let \(G^\ast\) be obtained from two disjoint copies of \(K_{3,3}\) by deleting edges \(x_1y_1,x_2y_2\) and adding \(x_1y_2,y_1x_2\). The resulting graph is cubic and bipartite, hence 2-connected, and \(\{x_1,y_1\}\) is a vertex cut.

Suppose that \(G^\ast\) admits a fall 3-colouring. Within each copy of \(K_{3,3}-xy\), the two non-terminal vertices in each part must have the same colour. The neighbourhood condition then forces the terminals \(x,y\) to have the third colour and their external neighbours to have different colours. This contradicts the equality of the terminal colours in the other copy. Therefore, Theorem 1.1 gives \(\chi_h^{(2)}(G^\ast)=4\).

The graph \(G^\ast\) is non-det-extremal: a perfect matching containing the two joining edges can be switched along a 4-cycle, reversing its sign. Thus, this example leaves open whether every det-extremal cubic bipartite graph of connectivity two admits a fall 3-colouring.

For cubic claw-free graphs, three colours suffice apart from $K_4$.
Indeed, let $G\ne K_4$ be such a graph. By Brooks' theorem~\cite{Brooks},
$G$ admits a proper 3-colouring. Each vertex has two adjacent neighbours,
which use the two colours different from its own. The colouring is
therefore 2-homogeneous, and $\hnum(G)=3$. This observation also occurs
in the proof of Li and Zhou's Theorem~3.1~\cite{LZ}. On the other hand,
$K_4$ admits no proper 2-homogeneous colouring, since the neighbours of
each vertex have three distinct colours in every proper colouring.

The general question of Lu\v{z}ar and Sot\'ak is the following.

\begin{question}[{\cite[Problem~5.1]{CC}}]
If a cubic graph admits a proper 2-homogeneous colouring, must
it admit one with at most four colours?
\end{question}

Theorem~\ref{main} answers this question for bipartite graphs. The
claw-free case has the stronger bound of three, whereas the question
remains open for general cubic graphs.

\end{document}